\documentclass[11pt]{amsart}

\usepackage[T1]{fontenc}
\usepackage[utf8]{inputenc}
\usepackage{amsmath,amssymb,mathtools}
\usepackage{microtype}
\usepackage[hmargin=1in, vmargin=1in]{geometry}
\usepackage[dvipsnames]{xcolor}
\usepackage[colorlinks=true,linkcolor=blue,citecolor=blue,urlcolor=blue]{hyperref}

\newcommand{\leanrepository}{https://github.com/pedrotradacete/FBLOrderContinuousLH}
\newcommand{\leanrevision}{167b82ad69724d2ddb2a701fc932e52d445faecc}
\DeclareRobustCommand{\leanref}[3]{%
  \nobreak\href{\leanrepository/blob/\leanrevision/FBLOrderContinuousLH/#1.lean\#L#2}{%
    \textcolor{RedOrange}{\textnormal{%
      {\scriptsize\ensuremath{\boldsymbol{\forall}}}\,\texttt{\detokenize{#3}}}}}%
  \space\ignorespaces}

\newtheorem{theorem}{Theorem}[section]
\newtheorem{proposition}[theorem]{Proposition}
\newtheorem{lemma}[theorem]{Lemma}
\newtheorem{corollary}[theorem]{Corollary}

\theoremstyle{definition}

\theoremstyle{remark}
\newtheorem{remark}[theorem]{Remark}

\numberwithin{equation}{section}

\newcommand{\FBL}[1]{\operatorname{FBL}[#1]}
\newcommand{\FBLp}[2]{\operatorname{FBL}^{(#1)}[#2]}
\newcommand{\FVL}[1]{\operatorname{FVL}[#1]}
\newcommand{\Hph}[1]{H[#1]}
\newcommand{\Cph}{C_{\mathrm{ph}}}
\newcommand{\bwstar}{bw^*}

\newcommand{\dist}{\operatorname{dist}}

\title[Order continuity in free Banach lattices]
{Order continuity of induced homomorphisms between free Banach lattices}

\author[Pariente-Mart\'in]{Ra\'ul Enrique Pariente-Mart\'in}
\address{Instituto de Ciencias Matem\'aticas (CSIC-UAM-UC3M-UCM)\\
Consejo Superior de Investigaciones Cient\'ificas\\
C/ Nicol\'as Cabrera, 13--15, Campus de Cantoblanco UAM\\
28049 Madrid, Spain.}
\email{ }

\author[Tradacete]{\\ Pedro Tradacete}
\address{Instituto de Ciencias Matem\'aticas (CSIC-UAM-UC3M-UCM)\\
Consejo Superior de Investigaciones Cient\'ificas\\
C/ Nicol\'as Cabrera, 13--15, Campus de Cantoblanco UAM\\
28049 Madrid, Spain.}
\email{pedro.tradacete@icmat.es}

\subjclass[2020]{Primary 46B42; Secondary 46A40, 47B60}

\keywords{Free Banach lattice, order-continuous operator, free vector lattice,
operator extension, bounded weak-star topology, $p$-compact sets, order-dense sublattice, majorizing sublattice}

\begin{document}

\begin{abstract}
Let $T\colon F\to E$ be a bounded operator between Banach spaces and let
$\overline T\colon \FBLp{p}{F}\to\FBLp{p}{E}$ be the lattice homomorphism
induced by $T$ between the corresponding free Banach lattices generated by $F$ and $E$.  We show that $\overline T$ is order continuous if and only if it is either zero or a lattice isomorphic embedding. For
$p=\infty$, this holds precisely when $T$ is zero or an isomorphic embedding, while for
$1\leq p<\infty$, apart from the zero operator,
$T$ must be an embedding and the pair $(T(F),E)$ must have the operator-extension property POE-$p$. In addition, we study when $\FBLp{p}{E}$ is order dense or majorizing in $\FBLp{q}{E}$ for $p<q$.
\end{abstract}

\maketitle

\section{Introduction}

The free Banach lattice generated by a Banach space was introduced in
\cite{ART18}, with its $p$-convex counterpart developed in
\cite{JardonSanchezEtAl2022}. These constructions extend the notion of free Banach lattice generated by a set given in \cite{dePagterWickstead2015}. More precisely, if $E$ is a Banach space and
$1\leq p\leq\infty$, then $\FBLp{p}{E}$ is a $p$-convex Banach lattice with
$p$-convexity constant one, together with a linear isometry
$\delta_E\colon E\to\FBLp{p}{E}$, satisfying the following
universal property: for every $p$-convex Banach lattice $X$ and every bounded
operator $S\colon E\to X$, there is a unique lattice homomorphism
$\hat S\colon\FBLp{p}{E}\to X$ such that
$\hat S\circ\delta_E=S$, and
$\|\hat S\|\leq M^{(p)}(X)\|S\|$, where $M^{(p)}(X)$ denotes the
$p$-convexity constant of $X$ (see \cite{JardonSanchezEtAl2022,OikhbergEtAl2026}). In particular, the assignment
$E\mapsto\FBLp{p}{E}$ is functorial: every bounded operator
$T\colon F\to E$ induces a unique lattice homomorphism
\[
 \overline T\colon \FBLp{p}{F}\longrightarrow \FBLp{p}{E},
\]
such that $\overline T\circ\delta_F=\delta_E\circ T,$ and
 with $\|\overline T\|=\|T\|$. This construction has motivated two closely related lines of research: understanding how geometric properties of $E$ are encoded by
$\FBLp{p}{E}$, and studying how properties of $T$ are inherited by
$\overline T$ (see, for example, \cite{Oikhberg2024,OikhbergEtAl2026}).  

Free Banach lattices have proved to be a useful tool in the understanding of the relation between Banach spaces and Banach lattices. 
They have been used to construct push-outs in Banach lattices, to obtain
amalgamation results and the construction of Banach lattices of universal
disposition and separably injective Banach lattices
\cite{AvilesTradacete2023}.  Recent developments also include the extension of
the theory to complex Banach spaces \cite{dHT}, Banach $f$-algebras \cite{MunozLahozTradacete2026}, non locally convex settings \cite{SalgueroTradaceteTrejo2025}, or the construction and
analysis of free and injective objects in the class of Banach lattices with
upper $p$-estimates \cite{GLTT}, among many others.

Several Banach space properties have been characterized in terms of the corresponding free Banach lattice. 
For example, finite dimensionality, separability, and lack of complemented copies of $\ell_1$ are all reflected respectively as existence of strong unit, quasi-interior points and order continuity of the dual of the corresponding free Banach lattice \cite[Section 9]{OikhbergEtAl2026}. Further properties, like reflexivity and approximation properties have also been characterized in these terms (see \cite{LaustsenTradacete2026,Oikhberg2024}).

Concerning operators, a number of properties are known to be preserved: $\overline T$ is
injective if and only if $T$ is injective, it has dense range if and only if
$T$ has dense range, and it is surjective if and only if $T$ is surjective
\cite[Proposition~3.2]{OikhbergEtAl2026}. The situation for isomorphic
embeddings is subtler. The operator $\overline T$ is an isomorphic lattice
embedding if and only if $T:F\rightarrow E$ is an isomorphic embedding and the pair
$(T(F),E)$ has the POE-$p$: there is a constant $C\geq1$ such that, for every
$n\in\mathbb N$, every operator $R\colon T(F)\to\ell_p^n$ extends to an
operator $\widetilde R\colon E\to\ell_p^n$ satisfying
$\|\widetilde R\|\leq C\|R\|$
\cite[Section~3]{OikhbergEtAl2026}.

In this context, it is natural to wonder under which conditions the lattice homomorphism $\overline T$ is order continuous. In \cite[Theorem 3.4]{OikhbergEtAl2026} it is shown that if $T:F\rightarrow E$ is an embedding, and $T(F)$ is complemented in $E$, then $\overline T$ is order continuous. In fact, as noticed by \cite{AzouziDhifaoui2025}, the original statement in \cite{OikhbergEtAl2026} containned the false claim that complementedness was not necessary (see also \cite{OikhbergEtAl2026Corrigendum}). Our aim in this paper is to give a complete characterization of order continuity of $\overline T$ in terms of properties of $T$ and the pair $(T(F),E)$, clarifying in particular the above issue.

Complementing earlier results in \cite{AzouziDhifaoui2025}, our first main result gives a complete answer for the case $p=\infty$. Throughout, an \emph{embedding} means an isomorphic embedding, or,
equivalently, an injective operator with closed range.

\begin{theorem}\label{thm:intro-infty}
\leanref{Order}{1804}{intro_infty}
Let $E$ and $F$ be Banach spaces and let $T\colon F\to E$ be bounded.  The
induced homomorphism
\[
 \overline T\colon \FBLp{\infty}{F}\longrightarrow
 \FBLp{\infty}{E}
\]
is order continuous if and only if $T$ is either zero or an embedding.
\end{theorem}

It was already shown in \cite{AzouziDhifaoui2025}, that if $T$ is an embedding, then $\overline T\colon \FBLp{\infty}{F}\longrightarrow \FBLp{\infty}{E}$ is order continuous. For proving the converse implication, here we exploit the key fact that a sequence $(u_k)\subset E$ without convergent subsequences satisfies $\bigwedge_{k=1}^n|\delta_{u_k}|\downarrow0$ in $\FBLp{\infty}{E}$ (see Proposition~\ref{prop:noncompact-minima}). In fact, for every finite $p$, this necessary condition remains valid:

\begin{proposition}\label{prop:necessary}
\leanref{Order}{1470}{necessary}
Let $E$ and $F$ be Banach spaces, let $1\leq p\leq\infty$, and let
$T\colon F\to E$ be bounded.  If
$\overline T\colon\FBLp{p}{F}\to\FBLp{p}{E}$ is order continuous, then $T$
is either zero or an embedding.
\end{proposition}

Although order continuity is defined using arbitrary nets, we note that in the present setting it is enough to test it with sequences; this useful preliminary is recorded in
Lemma~\ref{lem:sigma-order-continuity}. This fact together with an adaptation of the argument in \cite[Theorem~3.4]{OikhbergEtAl2026} allowed us to conclude that for the canonical inclusion $c_0\hookrightarrow\ell_\infty$, the induced homomorphism is order continuous for every $p$ (see Corollary~\ref{cor:c0}). This particular example shows that complementedness is not the necessary condition for order continuity. Nevertheless the pair $(c_0,\ell_\infty)$ is known to have POE-$p$. This fact, together with an earlier proof that for the span of the rademacher sequence in $L_1[0,1]$, the induced homomorphism is not order continuous, motivated us to look for a more general condition that would characterize order continuity of $\overline T$ for finite $p$.

\begin{theorem}\label{thm:intro-poe}
\leanref{Extension}{735}{intro_poe}
Let $1\leq p<\infty$ and let $T\colon F\to E$ be a bounded operator between
Banach spaces.  The induced homomorphism
\[
 \overline T\colon\FBLp{p}{F}\longrightarrow\FBLp{p}{E}
\]
is order continuous if and only if $T=0$, or $T$ is an embedding and
$(T(F),E)$ has POE-$p$.  Equivalently, $\overline T$ is order continuous if
and only if it is either zero or an isomorphic lattice embedding.
\end{theorem}

As mentioned above, \cite[Theorem~3.7]{OikhbergEtAl2026} identifies POE-$p$ with boundedness
below of the induced inclusion. Our main contribution is that this 
condition is actually equivalent to order continuity. An important ingredient of the proof is the topology induced by support functions of $p$-compact disks, which brings a new perspective for understanding the elements of the free $p$-convex Banach lattice. 

For $p<q$, we also consider the formal inclusion
\[
 J_E^{p,q}\colon\FBLp{p}{E}\longrightarrow\FBLp{q}{E},
\]
extending the identity on $E$. We show that the image of $J_E^{p,q}$ is a regular/order dense/majorizing sublattice only when $\FBLp{p}{E}$ is actually $q$-convex and hence $\FBLp{p}{E}$ and $\FBLp{q}{E}$ coincide up to equivalence of norms (see Theorem \ref{thm:cross-parameter-inclusion}).

The paper is organized as follows. Section~\ref{sec:preliminaries} recalls the basics and notation concerning free Banach lattices, establishes Lemma~\ref{lem:sigma-order-continuity} and some related reductions simplifying the test for order continuity, as well as the purely algebraic characterization of order continuity on the level of the free vector lattices (Theorem \ref{thm:fvl-characterization}). In Section~\ref{sec:induced} we provide the proofs of Theorem \ref{thm:intro-infty} and Proposition \ref{prop:necessary}.
Section~\ref{sec:poe} is devoted to the proof of Theorem~\ref{thm:intro-poe}. Finally, in Section~\ref{sec:canonical-inclusion} we analyze the order properties of $ J_E^{p,q}\colon\FBLp{p}{E}\longrightarrow\FBLp{q}{E}$.

\subsection*{Formalization} 
The theorems, propositions, lemmas, and corollaries in this paper have been
formalized in Lean~4.35.0 using \href{https://github.com/leanprover-community/mathlib4}{\texttt{Mathlib}} and the Banach lattice library \href{https://github.com/davidmunozlahoz/banlat}{\texttt{BanLat}}.
The code is pinned to Mathlib~4.35.0 (commit~\href{https://github.com/leanprover-community/mathlib4/commit/56a005855ea17e0b0958f73ec801ae140b0a1679}{\texttt{56a005855ea1}}) and BanLat commit~\href{https://github.com/davidmunozlahoz/banlat/commit/33c4b73fce845c4295614326193d4e52111be762}{\texttt{33c4b73fce84}}.
The formalization is available in the \href{\leanrepository}{\texttt{FBLOrderContinuousLH}} repository.
Each statement is accompanied by a linked $\forall$ symbol and its Lean
declaration name. These links point to the corresponding source line at a
fixed revision of the repository, so that they remain stable as the code evolves.

\subsection*{Declaration on use of AI}
GPT-6 Astra was used for providing a first version of the proof of Theorem \ref{thm:intro-poe}. This proof was digested and subsequently simplified by the authors. GPT-6 Astra and GPT-6 Sol were also used for the Lean formalization of the results in the paper.

\section{Preliminaries and criteria for order continuity}\label{sec:preliminaries}

All Banach spaces and vector lattices are real.  If $(x_\alpha)$ is a net in
a vector lattice, we write $x_\alpha\downarrow0$ when it is decreasing and
has infimum zero.  A positive linear operator $S\colon X\to Y$ is
\emph{order continuous} if $x_\alpha\downarrow0$ implies
$Sx_\alpha\downarrow0$, and it is \emph{$\sigma$-order continuous} if the
same condition is required only for sequences, see
\cite[Chapter~1]{AliprantisBurkinshaw2006} for background.

Let $E$ be a Banach space.  Let $\Hph{E}$ denote the vector lattice of all
positively homogeneous functions $f\colon E^*\to\mathbb R$, endowed with
pointwise order. For $x\in E$, let us consider the evaluation functional $\delta_x\in\Hph{E}$ defined by
\[
 \delta_x(x^*)=x^*(x),
\]
for $x^*\in E^*$. The sublattice of $\Hph{E}$ generated by $\{\delta_x:x\in E\}$ coincides with the free vector lattice generated by $E$, $\FVL{E}$, see \cite{ART18,baker}.  

For $1\leq p<\infty$ and $f\in\Hph{E}$, let
\begin{equation}\label{eq:fblp-norm}
 \|f\|_{\FBLp{p}{E}}
 =\sup\left\{
 \left(\sum_{k=1}^n|f(x_k^*)|^p\right)^{1/p}:
 \sup_{x\in B_E}\sum_{k=1}^n|x_k^*(x)|^p\leq1
 \right\},
\end{equation}
and let 
\[
H_p[E]=\{f\in \Hph{E}:\|f\|_{\FBLp{p}{E}}<\infty\}.
\] 
This is a Banach lattice with the above norm and pointwise order. The free Banach lattice $\FBLp{p}{E}$ can be represented as the completion of $\FVL{E}$ with respect to this norm, in other words it corresponds to the closed sublattice generated by $\{\delta_x:x\in E\}$ in $H_p[E]$ (see \cite[Theorem~2.1]{OikhbergEtAl2026}).  The map
$\delta_E\colon x\mapsto\delta_x$ is a linear isometry.  We also write
$\FBL{E}=\FBLp{1}{E}$.

For $p=\infty$, the free $\infty$-convex Banach lattice can be identified with the Banach lattice of continuous positively homogeneous functions on the closed unit ball of the dual space, endowed with the weak-star topology, which we denote as
\begin{equation}\label{eq:infty-representation}
 \FBLp{\infty}{E}=\Cph(B_{E^*}),
 \qquad
 \|f\|_{\FBLp{\infty}{E}}=\sup_{x^*\in B_{E^*}}|f(x^*)|,
\end{equation}
\cite[Proposition~2.2]{OikhbergEtAl2026}.  For every $p$, the elements of
$\FBLp{p}{E}$ are weak-star continuous on bounded subsets of $E^*$ and hence
continuous for the bounded weak-star topology. We will use the canonical contractive
injective lattice homomorphisms
\[
 J_E^{p,q}\colon\FBLp{p}{E}\longrightarrow\FBLp{q}{E}
 \qquad(1\leq p<q\leq\infty),
\]
(see \cite[Proposition~7.3]{JardonSanchezEtAl2022}).

The universal property gives, for each bounded operator between Banach spaces $T\colon F\to E$, a unique
lattice homomorphism $\overline T$ satisfying
$\overline T\delta_x=\delta_{Tx}$.  In the functional representation given above, this corresponds to the
composition operator
\begin{equation}\label{eq:composition}
 (\overline T f)(x^*)=f(T^*x^*)
\end{equation}
for $f\in\FBLp{p}{F}$ and $x^*\in E^*$.

We shall also make repeated use of the elementary lattice inequalities
\begin{equation}\label{eq:lattice-inequalities}
 \bigl||a|-|b|\bigr|\leq|a-b|,
 \qquad
 0\leq a-a\wedge b=(a-b)^+\leq|a-b|.
\end{equation}





Recall that a vector lattice has the \emph{countable supremum property} if every subset having a supremum contains a countable subset with the same supremum.

\begin{lemma}\label{lem:csp}
\leanref{Order}{453}{csp}
For every Banach space $E$ and every $1\leq p\leq\infty$, the lattice
$\FBLp{p}{E}$ has the countable supremum property.
\end{lemma}

\begin{proof}
By \cite{AvilesPlebanekRodriguez2018}, the Banach lattice $\Cph(B_{E^*})$ satisfies the $\sigma$-bounded chain condition. Therefore, so does every (not necessarily closed) sublattice, and in particular $\FBLp{p}{E}$. Thus every disjoint family of positive elements must be countable. The countable supremum property now follows from \cite[Theorem~29.3]{LuxemburgZaanen1971}.
\end{proof}

The countable supremum property allows us to reduce the order continuity test to sequences, as recorded in the next proposition.

\begin{proposition}\label{prop:csp-operator}
\leanref{Order}{92}{csp_operator}
Let $X$ be a vector lattice with the countable supremum property, let $Y$ be
a vector lattice, and let $S\colon X\to Y$ be a positive operator.  Then $S$ is order
continuous if and only if it is $\sigma$-order continuous.
\end{proposition}

\begin{proof}
We only prove the non-trivial implication. Suppose that $S$ is
$\sigma$-order continuous and let $x_\alpha\downarrow0$ in $X$.  Fix
$\alpha_0$.  Since
\[
 \sup_{\alpha\geq\alpha_0}(x_{\alpha_0}-x_\alpha)=x_{\alpha_0},
\]
the countable supremum property yields indices $\alpha_n\geq\alpha_0$ such
that $\inf_n x_{\alpha_n}=0$.  As the original index set is directed, we may
choose successively $\beta_n$ so that
\[
 \beta_n\geq\beta_{n-1},\alpha_n.
\]
Then $x_{\beta_n}\downarrow0$.  Hence
$Sx_{\beta_n}\downarrow0$.  Any positive lower bound of the decreasing net
$(Sx_\alpha)$ is a lower bound of this sequence and is therefore zero.  Thus
$Sx_\alpha\downarrow0$, as required.
\end{proof}

For convenience and further reference we state next the following particular case, which will be used in the sequel.

\begin{lemma}\label{lem:sigma-order-continuity}
\leanref{Order}{494}{sigma_order_continuity}
Let $E$ be a Banach space, $1\leq p\leq\infty$, and let $Y$ be a Banach
lattice.  A positive linear operator
\[
 S\colon \FBLp{p}{E}\longrightarrow Y
\]
is order continuous if and only if it is $\sigma$-order continuous.
\end{lemma}

\begin{proof}
Combine Lemma~\ref{lem:csp} and Proposition~\ref{prop:csp-operator}.
\end{proof}

The sequential test can be restricted further to the norm-dense vector
sublattice $\FVL{F}$ as follows.

\begin{proposition}\label{prop:fvl-sequential-test}
\leanref{Order}{600}{fvl_sequential_test}
Let $1\leq p\leq\infty$, and let $T\colon F\to E$ be a bounded operator
between Banach spaces. The induced lattice homomorphism
$\overline T\colon\FBLp{p}{F}\to\FBLp{p}{E}$ is order continuous if and
only if, for every sequence $(f_n)\subseteq\FVL{F}_+$ with
$f_n\downarrow0$ in $\FBLp{p}{F}$, one has
$\overline T f_n\downarrow0$ in $\FBLp{p}{E}$.
\end{proposition}

\begin{proof}
Only the reverse implication requires proof. By
Lemma~\ref{lem:sigma-order-continuity}, it suffices to prove
$\sigma$-order continuity. Suppose otherwise. Then there are
$h_n\downarrow0$ in $\FBLp{p}{F}$ and
$0<u\in\FBLp{p}{E}$ such that
$u\leq\overline T h_n$ for every $n$.
Taking positive parts of norm approximants from $\FVL{F}$, choose
$g_n\in\FVL{F}_+$ such that
\[
 \sum_{n=1}^{\infty}\|g_n-h_n\|_{\FBLp{p}{F}}
 <\frac{\|u\|_{\FBLp{p}{E}}}{1+\|\overline T\|}.
\]
Put $f_n=\bigwedge_{k=1}^n g_k$. Then $(f_n)$ is a decreasing
sequence in $\FVL{F}_+$. If $0\leq h\leq f_n$ for every $n$, then,
for each fixed $m$ and every $n\geq m$,
\[
 h\leq g_n\leq h_m+|g_n-h_n|.
\]
The positive elements $h_m+|g_n-h_n|-h$ converge in norm to
$h_m-h$. Since the positive cone is norm closed, $h\leq h_m$.
As this holds for every $m$, we have $h=0$. Thus
$f_n\downarrow0$ in $\FBLp{p}{F}$.

The positive series
\[
 e=\sum_{n=1}^{\infty}|\overline T(g_n-h_n)|
\]
converges in $\FBLp{p}{E}$ and satisfies
\[\|e\|_{\FBLp{p}{E}}\leq\|\overline T\|\sum_n\|g_n-h_n\|_{\FBLp{p}{F}}<\|u\|_{\FBLp{p}{E}}.\]
Consequently $v=(u-e)^+>0$, since otherwise $u\leq e$ and hence
$\|u\|_{\FBLp{p}{E}}\leq\|e\|_{\FBLp{p}{E}}$. For every $n$,
\[
 u-e\leq\overline T h_n-|\overline T(g_n-h_n)|
       \leq\overline T g_n.
\]
As $\overline T g_n\geq0$, it follows that $v\leq\overline T g_n$.
Hence
\[
 0<v\leq\bigwedge_{k=1}^n\overline T g_k
       =\overline T f_n,
\]
for every $(n\in\mathbb N$, contrary to the assumed test. This proves the assertion.
\end{proof}

Recall that a closed subspace $Y$ of a Banach space $X$ is a
\emph{Banach-space ideal} in $X$ if its annihilator
\[
 Y^\perp=\{x^*\in X^*:x^*|_Y=0\}
\]
is the kernel of a contractive projection on $X^*$.  Equivalently, there
is a linear isometry $\Phi\colon Y^*\to X^*$ such that
$(\Phi y^*)|_Y=y^*$ for every $y^*\in Y^*$; thus $\Phi$ selects
Hahn--Banach extensions linearly.  This notion involves no lattice
structure and should be distinguished from that of an order ideal;
see \cite[Section~4.2]{OikhbergEtAl2026}.

The following reduction shows that failure of order continuity is always
witnessed inside separable Banach-space ideals.

\begin{lemma}\label{lem:separable-reduction}
\leanref{Order}{1073}{separable_reduction}
Let $T\colon F\to E$ be a bounded operator and let $1\leq p\leq\infty$.
If
\[
 \overline T\colon\FBLp{p}{F}\longrightarrow\FBLp{p}{E}
\]
is not order continuous, then there are separable Banach-space ideals
$F_0\subseteq F$ and $E_0\subseteq E$, with $T(F_0)\subseteq E_0$, such
that the homomorphism induced by the restriction
\[
 T_0=T|_{F_0}\colon F_0\longrightarrow E_0
\]
is not order continuous.  If $F$ is separable, one may take $F_0=F$.
\end{lemma}

\begin{proof}
By Lemma~\ref{lem:sigma-order-continuity}, there are
$f_n\downarrow0$ in $\FBLp{p}{F}$ and $0<g\in\FBLp{p}{E}$ such that for $n\in\mathbb N$
\begin{equation}\label{eq:separable-reduction-witness}
 g\leq\overline T f_n.
\end{equation}
For each $n$, choose a sequence in $\FVL{F}$ converging to $f_n$, and let
$A_F\subseteq F$ be the countable set of all generators occurring in these
approximations.  Choose a separable Banach-space ideal $F_0$ of $F$
containing $A_F$ \cite{HeinrichMankiewicz1982}.  When $F$ is separable, take
$F_0=F$.

Similarly, choose a countable set $A_E\subseteq E$ containing all the
generators occurring in a sequence of lattice combinations that approximate
$g$, and let $E_0$ be a separable Banach-space ideal of $E$ containing
$A_E\cup T(F_0)$.  Let $\iota_F\colon F_0\hookrightarrow F$ and
$\iota_E\colon E_0\hookrightarrow E$ denote the inclusions.  The induced
lattice homomorphisms
\[
 \overline{\iota_F}\colon\FBLp{p}{F_0}\longrightarrow\FBLp{p}{F},
 \qquad
 \overline{\iota_E}\colon\FBLp{p}{E_0}\longrightarrow\FBLp{p}{E}
\]
are isometric lattice embeddings
\cite[Theorem~3.7 and Corollary~4.14]{OikhbergEtAl2026}.  Their ranges are
therefore closed.  By construction, there are
$\widetilde f_n\in\FBLp{p}{F_0}$ and
$\widetilde g\in\FBLp{p}{E_0}$ such that
\[
 \overline{\iota_F}\widetilde f_n=f_n,
 \qquad \overline{\iota_E}\widetilde g=g.
\]

We have $\widetilde f_n\downarrow0$ in $\FBLp{p}{F_0}$.  Indeed, any
positive common lower bound would be carried by $\overline{\iota_F}$ to a
positive common lower bound of $(f_n)$ and would consequently vanish.
Moreover, uniqueness
in the universal property gives the commutative identity
\[
 \overline{\iota_E}\overline{T_0}=\overline T\overline{\iota_F}.
\]
Since an injective lattice homomorphism reflects order,
\eqref{eq:separable-reduction-witness} yields
\[
 0<\widetilde g\leq\overline{T_0}\widetilde f_n
 \qquad(n\in\mathbb N).
\]
Thus $\overline{T_0}$ is not $\sigma$-order continuous and hence, by
Lemma~\ref{lem:sigma-order-continuity}, is not order continuous.
\end{proof}

For a vector space $E$, let $E^\#$ denote its algebraic dual.  The free
vector lattice $\FVL{E}$ can be realized as the sublattice of
$\mathbb R^{E^\#}$ generated by the evaluations $\delta_x$.  If
$T\colon F\to E$ is linear, its extension is
\begin{equation}\label{eq:fvl-composition}
 \overline T^{\mathrm{FVL}}f=f\circ T^\#,
\end{equation}
where $T^\#\colon E^\#\to F^\#$ is the algebraic adjoint.

The proof from \cite[Theorem~3.4]{OikhbergEtAl2026} that if $F$ is a complemented subspace of $E$, then the induced homomorphism $\overline\iota\colon\FBLp{p}{F}\to\FBLp{p}{E}$ is order continuous can be straightforwardly extended to the free vector lattice setting:

\begin{lemma}\label{lem:regular-inclusion}
\leanref{Order}{1964}{regular_inclusion}
Let $F$ be a vector subspace of $E$ and let
$\iota\colon F\hookrightarrow E$ be the inclusion. The extension
$\overline\iota^{\mathrm{FVL}}\colon\FVL{F}\to\FVL{E}$ is order
continuous.
\end{lemma}

As a consequence of this, we obtain the following characterization of order continuity for extensions of linear operators to free vector lattices.

\begin{theorem}\label{thm:fvl-characterization}
\leanref{Order}{1971}{fvl_characterization}
Let $T\colon F\to E$ be a nonzero linear operator between vector spaces.
Then
\[
 \overline T^{\mathrm{FVL}}\colon\FVL{F}\longrightarrow\FVL{E}
\]
is order continuous if and only if $T$ is injective.
\end{theorem}

\begin{proof}
Suppose first that $T$ is injective.  Factor it as
$T=\iota\circ T_0$, where $T_0\colon F\to T(F)$ is a linear isomorphism and
$\iota\colon T(F)\hookrightarrow E$ is the inclusion.  The extension of
$T_0$ is a lattice isomorphism and the extension of $\iota$ is order
continuous by Lemma~\ref{lem:regular-inclusion}.  Hence, the composition
$\overline T^{\mathrm{FVL}}$ is order continuous.

Conversely, suppose that $T$ is not injective.  Since $T\neq0$, there are
linearly independent $v_1,v_2\in F$ such that
$Tv_1=Tv_2=u\neq0$.  For $n\geq1$ set
\[
 f_n=\bigl(|\delta_{v_1}+\delta_{v_2}|
              -n|\delta_{v_1}-\delta_{v_2}|\bigr)^+
 \in\FVL{F}.
\]
The sequence is decreasing.  If $0\leq f\leq f_n$ for every $n$, then $f$
vanishes at every $x^\#\in F^\#$ satisfying
$x^\#(v_1-v_2)\neq0$.  This set is dense for the topology of pointwise
convergence on $F$, while every member of $\FVL{F}$ is continuous for that
topology.  Hence $f=0$, and therefore $f_n\downarrow0$.  However,
\[
 \overline T^{\mathrm{FVL}}f_n
 =\bigl(2|\delta_u|-n|\delta_u-\delta_u|\bigr)^+
 =2|\delta_u|
\]
for every $n$.  Thus $\overline T^{\mathrm{FVL}}$ is not order continuous.
\end{proof}

Note that the zero lattice homomorphism is order continuous even though its underlying linear map, $T=0$, is not injective.

\section{Order continuity of induced homomorphisms only holds for embeddings}\label{sec:induced}

We first note that order continuity of $\overline{T}:\FBLp{p}{F}\rightarrow \FBLp{p}{E}$ detects injectivity of the underlying operator for every $1\leq p\leq\infty$.

\begin{proposition}\label{prop:injectivity}
\leanref{Order}{1270}{injectivity}
Let $T\colon F\to E$ be a nonzero bounded operator and let
$1\leq p\leq\infty$.  If $\overline T\colon\FBLp{p}{F}\to\FBLp{p}{E}$ is
order continuous, then $T$ is injective.
\end{proposition}

\begin{proof}
If $T$ is not injective, choose linearly independent $v_1,v_2\in F$ with
$Tv_1=Tv_2=u\neq0$ and define $f_n$ as in the proof of
Theorem~\ref{thm:fvl-characterization}.  Every member of $\FBLp{p}{F}$ is
weak-star continuous on bounded subsets of $F^*$ \cite[Section 2]{OikhbergEtAl2026}.  Therefore the same
density argument, now on each dual ball, shows that $f_n\downarrow0$ in
$\FBLp{p}{F}$.  However, $\overline T f_n=2|\delta_u|$ for all $n$, contradicting
order continuity.
\end{proof}

We next isolate the property of sequences in the underlying space used to detect that the range is not closed.

\begin{proposition}\label{prop:noncompact-minima}
\leanref{Order}{1408}{noncompact_minima}
Let $(u_n)$ be a bounded sequence in a Banach space $F$ with no
norm-convergent subsequence.  Let
\[
 f_k=\bigwedge_{n=1}^k|\delta_{u_n}|.
\]
For every $1\leq p\leq\infty$, we have that $ f_k\downarrow0$ in $\FBLp{p}{F}.$
\end{proposition}

\begin{proof}
There is $\varepsilon_0>0$ such that for every finite-dimensional subspace
$M\subseteq F$ some $u_k$ satisfies
$\dist(u_k,M)>\varepsilon_0$.  Indeed, failure of this would make
$\{u_n:n\in\mathbb N\}$ totally bounded and hence relatively compact.

Let $C=\sup_n\|u_n\|$ and
$\lambda=\varepsilon_0/(\varepsilon_0+C)$.  For each $k$ set
\[
 A_k=\{x^*\in B_{F^*}:x^*(u_k)=0\}.
\]
We claim that $\lambda B_{F^*}$ lies in the weak-star closure of
$\bigcup_kA_k$.  Fix $x^*\in B_{F^*}$ and a basic weak-star neighbourhood of
$\lambda x^*$ determined by $z_1,\ldots,z_m\in F$.  Put
$M=\operatorname{span}\{z_1,\ldots,z_m\}$ and choose $k$ with
$\dist(u_k,M)>\varepsilon_0$.  Hahn--Banach gives $w^*\in F^*$ such that
\[
 w^*|_M=0,\qquad
 w^*(u_k)=-\lambda x^*(u_k),\qquad
 \|w^*\|\leq\frac{\lambda C}{\varepsilon_0}.
\]
Then $y^*=\lambda x^*+w^*$ belongs to $B_{F^*}$, agrees with $\lambda x^*$
on $M$, and satisfies $y^*(u_k)=0$.  This proves the claim.

Now, let $f_k=\bigwedge_{n=1}^k|\delta_{u_n}|$. Clearly, $(f_k)$ is a decreasing sequence. If $0\leq g\leq f_k$ for all $k$, then $g$ vanishes on every $A_k$.  By
weak-star continuity on $B_{F^*}$ it vanishes on $\lambda B_{F^*}$, and
positive homogeneity gives $g=0$.  Hence $f_k\downarrow0$.
\end{proof}

Now, we can prove the main results of this section.

\begin{proof}[Proof of Proposition~\ref{prop:necessary}]
If $T=0$, there is nothing to prove.  Suppose that $T\neq0$.  By
Proposition~\ref{prop:injectivity}, $T$ is injective.  Suppose that $T$ is
not bounded below.  Choose $x_n\in S_F$ with $Tx_n\to0$, and choose
$x_0\in F$ with $Tx_0\neq0$.  No subsequence of $(x_n)$ converges in norm:
otherwise its limit would belong to $S_F\cap\ker T$.  Thus
$u_n=x_0+x_n$ has no norm-convergent subsequence.

Choose $e^*\in E^*$ with $e^*(Tx_0)=1$.  Passing to a subsequence, we may
assume that
\begin{equation}\label{eq:subsequence-control}
 |e^*(Tu_n)|\geq\frac12,
 \qquad
 \|Tx_{n+1}-Tx_n\|\leq2^{-n}
 \quad(n\geq1).
\end{equation}
Put $f_k=\bigwedge_{n=1}^k|\delta_{u_n}|$.  By
Proposition~\ref{prop:noncompact-minima}, $f_k\downarrow0$ in
$\FBLp{\infty}{F}$.  Since $\FBLp{p}{F}$ is a (in general, non-closed) sublattice of
$\FBLp{\infty}{F}$, the same sequence decreases to zero in
$\FBLp{p}{F}$.

Set $g_k=\overline T f_k$.  Since $\overline T$ preserves lattice
operations,
\[
 g_n=\bigwedge_{k=1}^n|\delta_{Tu_k}|\leq|\delta_{Tu_n}|,
 \qquad g_{n+1}=g_n\wedge|\delta_{Tu_{n+1}}|.
\]
The identity $x-x\wedge y=(x-y)^+$ and monotonicity of the positive
part, together with \eqref{eq:lattice-inequalities}, give
\begin{align*}
 0\leq g_n-g_{n+1}
 &=(g_n-|\delta_{Tu_{n+1}}|)^+ \\
 &\leq(|\delta_{Tu_n}|-|\delta_{Tu_{n+1}}|)^+ \\
 &\leq\bigl||\delta_{Tu_n}|-|\delta_{Tu_{n+1}}|\bigr| \\
 &\leq|\delta_{T(u_n-u_{n+1})}|,
\end{align*}
and hence, by \eqref{eq:subsequence-control},
\[
 \|g_n-g_{n+1}\|\leq\|Tx_n-Tx_{n+1}\|\leq2^{-n}.
\]
It follows that $(g_n)$ converges in norm to some $g\geq0$.  Since it is
decreasing, $g\leq g_n$ for all $n$.  On the other hand,
\[
 g_n(e^*)=\min_{1\leq k\leq n}|e^*(Tu_k)|\geq\frac12,
\]
so $g(e^*)\geq1/2$ and $g\neq0$.  Thus $(g_n)$ does not decrease to zero,
contrary to the order continuity of $\overline T$.  Therefore $T$ is bounded
below and hence an embedding.
\end{proof}

\begin{proof}[Proof of Theorem~\ref{thm:intro-infty}]
The necessity follows from Proposition~\ref{prop:necessary}.  Conversely,
the conclusion is immediate if $T=0$. If $T$ is an embedding, then $\overline T$ is order continuous by \cite{AzouziDhifaoui2025}.
\end{proof}

As a first consequence, we obtain an embedding whose extension is order
continuous for every $1\leq p\leq\infty$, despite its range is uncomplemented in the ambient space.

\begin{corollary}\label{cor:c0}
\leanref{Extension}{1312}{c0}
Let $\iota\colon c_0\hookrightarrow\ell_\infty$ be the canonical inclusion.
For every $1\leq p\leq\infty$, the induced lattice homomorphism
\[
 \overline\iota\colon\FBLp{p}{c_0}\longrightarrow
 \FBLp{p}{\ell_\infty}
\]
is order continuous.
\end{corollary}

\begin{proof}
Fix $1\leq p\leq\infty$ and suppose that $\overline \iota$ is not order continuous.  Lemma~\ref{lem:separable-reduction}, with $F_0=F=c_0$, provides a separable Banach-space ideal $E_0$ of $\ell_\infty$ containing $\iota(F)$ such that the extension of $\iota\colon c_0\to E_0$ is not order continuous.  Put $M=\iota(c_0)$. By Sobczyk's theorem \cite[Theorem~2.5.8]{AlbiacKalton2016} $M$ is complemented in $E_0$.  Factor $\iota\colon c_0\to E_0$ through $M$.  The extension of the isomorphism $c_0\to M$ is a lattice isomorphism, while the extension of the inclusion $M\hookrightarrow E_0$ is order continuous by \cite[Theorem 3.4]{OikhbergEtAl2026}.  Their composition is therefore order continuous, a contradiction.
\end{proof}

Actually, the argument of Corollary~\ref{cor:c0} works for any separable Banach space $F$ that is separably injective, and any embedding $T\colon F\to E$ into a Banach space $E$. Although this corollary could also be later obtained as a particular case of Theorem \ref{thm:intro-poe}, we have decided to include its details for motivation and exposition purposes.

\section{Operator extension property and order continuity of induced homomorphisms}
\label{sec:poe}

Throughout this section, $F$ is a closed subspace of a Banach space $E$,
and $\iota\colon F\hookrightarrow E$ is the inclusion. Unless explicitly
stated otherwise, we assume that $1\leq p<\infty$ and write $p'$ for the
conjugate exponent, with $p'=\infty$ when $p=1$. Hahn--Banach makes
$\iota^*$ surjective, so
$\overline\iota\colon\FBLp{p}{F}\to\FBLp{p}{E}$ is always
an injective lattice homomorphism. Its range need not be norm closed. Recall
that a sublattice is \emph{regular} if its inclusion is order continuous.
We shall prove that, for these induced inclusions, regularity is equivalent
to closedness of the range.

\begin{theorem}\label{thm:poe-inclusion}
\leanref{Extension}{537}{poe_inclusion}
Let $1\leq p<\infty$ and let $F$ be a closed subspace of a Banach space
$E$.  The following are equivalent:
\begin{enumerate}
\renewcommand{\labelenumi}{(\roman{enumi})}
\item $\overline\iota\colon\FBLp{p}{F}\to\FBLp{p}{E}$ is order continuous;
\item $\overline\iota$ is $\sigma$-order continuous;
\item $\overline\iota$ is bounded below;
\item $(F,E)$ has POE-$p$.
\end{enumerate}
If these conditions fail, there are $f_N\in\FVL{F}_+$ and
$0<G\in\FBLp{p}{E}$ such that
\[
 f_N\downarrow0\quad\text{in }\FBLp{p}{F},
 \qquad
 \overline\iota f_N\longrightarrow G
 \quad\text{in the norm of }\FBLp{p}{E}.
\]
\end{theorem}

The equivalence of (iii) and (iv) is the operator-extension theorem
\cite[Theorem~3.7 and Definition~3.8]{OikhbergEtAl2026}; the equivalence of
(i) and (ii) is also contained in Lemma~\ref{lem:sigma-order-continuity}.
We will next prove (iii)$\Rightarrow$(i) and (ii)$\Rightarrow$(iii). 

\subsection{Control with support functions}

For a nonempty bounded set $A$ in a Banach space $E$, write its support function as
\[
 h_A(x^*)=\sup_{x\in A}|x^*(x)|
\]
for $x^*\in E^*.$ This is a norm-continuous seminorm on $E^*$.   

Our first lemma provides an equivalent form of the free $p$-convex norm.
For $p=1$ it can be related to the maximal lattice norm given in
\cite[Theorem~3.1]{Troitsky2019}. For general $p$, compare also with the unit-ball
description in \cite[Proposition~2.4]{Oikhberg2024}.  We give a direct proof below. 


\begin{lemma}\label{lem:majorant-norm}
\leanref{Support}{823}{majorant_norm}
For every Banach space $E$, every $1\leq p<\infty$, and every
$v\in\FVL{E}$,
\begin{equation}\label{eq:majorant-norm}
 \|v\|_{\FBLp{p}{E}}
 =\inf\left\{\left(\sum_{j=1}^m\|x_j\|^p\right)^{1/p}:
 x_j\in E,\quad
 |v|\leq\left(\sum_{j=1}^m|\delta_{x_j}|^p\right)^{1/p}\right\}.
\end{equation}
\end{lemma}

\begin{proof}
Let $\mathcal A_E$ be the vector lattice of functions
$H(\delta_{x_1},\ldots,\delta_{x_m})$, where $m<\infty$, $x_j\in E$, and
$H\colon\mathbb R^m\to\mathbb R$ is continuous and positively homogeneous.
The usual homogeneous functional calculus realizes $\mathcal A_E$ inside
$\FBLp{p}{E}$.  Define $\nu_p$ on $\mathcal A_E$ by the right-hand side of
\eqref{eq:majorant-norm}.  Boundedness of $H$ on the unit sphere shows that
every such function has a finite majorant of the indicated form.
Homogeneity and monotonicity are immediate, and
$\|v\|_{\FBLp{\infty}{E}}\leq\nu_p(v)$ and
$\nu_p(\delta_x)=\|x\|$ follow from Hahn--Banach.

To check the triangle inequality, for $p=1$ is immediate. When $p>1$, suppose that $|f|\leq \left(\sum_{j=1}^m|\delta_{x_j}|^p\right)^{1/p}=a$ and
$|g|\leq \left(\sum_{j=1}^m|\delta_{y_j}|^p\right)^{1/p}=b$ with $A=\left(\sum_{j=1}^m\|x_j\|^p\right)^{1/p}$, $B=\left(\sum_{j=1}^m\|y_j\|^p\right)^{1/p}$, and assume $A,B>0$.  Weighted
H\"older gives
\[
 a+b\leq\left(t^{1-p}a^p+(1-t)^{1-p}b^p\right)^{1/p}
 \qquad(0<t<1).
\]
The expression on the right is again a finite majorant.  With
$t=A/(A+B)$, its corresponding contribution is $A+B$ and taking infima proves subadditivity. Thus $\nu_p$ is a lattice norm. 

Moreover, concatenating
majorants of $f_1,\ldots,f_s$ gives
\[
 \nu_p\left(\left(\sum_{k=1}^s|f_k|^p\right)^{1/p}\right)
 \leq\left(\sum_{k=1}^s\nu_p(f_k)^p\right)^{1/p}.
\]
The pointwise $p$-sum operation agrees with the homogeneous functional
calculus for this norm: approximate its defining function uniformly by
homogeneous lattice combinations, with the error bounded by a multiple of
$\sum_k|f_k|$.  The same approximation argument shows that
$\FVL{E}$ is $\nu_p$-dense in $\mathcal A_E$.

It follows that the completion $Z$ of $(\mathcal A_E,\nu_p)$ is a
$p$-convex Banach lattice with constant $1$.  The isometry
$x\mapsto\delta_x\in Z$ extends to a contractive lattice homomorphism
$\FBLp{p}{E}\to Z$.  On $\FVL{E}$ this extension is the identity, whence
$\nu_p(v)\leq\|v\|_{\FBLp{p}{E}}$ for $v\in \FVL{E}$.  Conversely, $p$-convexity in $\FBLp{p}{E}$
gives $\|v\|_{\FBLp{p}{E}}\leq(\sum_j\|x_j\|^p)^{1/p}$ for every admissible
majorant.  This proves \eqref{eq:majorant-norm}.
\end{proof}

In the following sections, for $1\leq p<\infty$, a \emph{$p$-compact disk}
in $E$ means a set of the form
\begin{equation}\label{eq:summable-disk}
 K=\left\{\sum_{j=1}^{\infty}t_jx_j:
 (t_j)\in B_{\ell_{p'}}\right\},
 \qquad \text{for some } (x_j)\subset E \text{ such that }\sum_{j=1}^{\infty}\|x_j\|^p<\infty.
\end{equation}
Such a set is absolutely convex and norm compact. Indeed,
$B_{\ell_{p'}}$ is compact for coordinatewise convergence, including when
$p'=\infty$, and H\"older bounds the norm of the series tail uniformly by
$(\sum_{j>N}\|x_j\|^p)^{1/p}$.  Its support function satisfies
\begin{equation}\label{eq:summable-support}
 h_K=\left(\sum_{j=1}^{\infty}|\delta_{x_j}|^p\right)^{1/p}
 \in\FBLp{p}{E},
 \qquad \|h_K\|_{\FBLp{p}{E}}\leq\left(\sum_j\|x_j\|^p\right)^{1/p}.
\end{equation}
The displayed $p$-sum is a norm limit of its finite partial $p$-sums.
Combining finitely many such support functions by taking a $p$-sum,
and multiplying them by positive scalars, again gives support functions
of $p$-compact disks. For notational convenience, an \emph{$\infty$-compact disk} in $E$ means a
nonempty absolutely convex norm-compact subset of $E$. Let us recall that the theory of $p$-compact sets arose as a natural generalization of compact sets in \cite{SinhaKarn2002}.

For $1\leq p\leq\infty$, let $\tau_p$ denote the locally convex topology on $E^*$
of uniform convergence on $p$-compact disks; equivalently, $\tau_p$ is the topology
generated by the seminorms $h_K$, with $K$ ranging over the $p$-compact
disks in $E$.

Observe that, for every nonempty norm-compact set $K\subset E$,
the closed absolutely convex hull 
$K^a=\overline{\operatorname{absconv}}^{\|\cdot\|}K$ is norm compact and
\[
 h_K=h_{K^a}.
\]
Hence the bounded weak-star topology $\bwstar$ on $E^*$, which is the
topology of uniform convergence on norm-compact subsets of $E$, is
equivalently generated by the seminorms $h_K$ with $K$ ranging over the
$\infty$-compact disks. Therefore
$
\tau_\infty=\bwstar.
$

For $1\leq p<\infty$, every $p$-compact disk is norm compact, whereas for
each $x\in E$ the disk $\{tx:|t|\leq1\}$ is $p$-compact and has support
function $|\delta_x|$. Hence
$$
\sigma(E^*,E)\subseteq\tau_p\subseteq\tau_\infty=\bwstar.
$$
In particular, each evaluation $\delta_x$ is $\tau_p$-continuous, and
therefore every element of $\FVL{E}$ is also $\tau_p$-continuous.

The following approximation argument allows us to extend $\tau_p$-continuity to elements in $\FBLp{p}{E}$.

\begin{lemma}\label{lem:summable-local-control}
\leanref{Support}{1858}{summable_local_control}
Let $1\leq p\leq\infty$, $g\in\FBLp{p}{E}$, $x_0^*\in E^*$ and
$\varepsilon>0$. There is a $p$-compact disk $K\subset E$ such that
$$
 h_K(x^*-x_0^*)<1
 \quad\Longrightarrow\quad
 |g(x^*)-g(x_0^*)|<\varepsilon.
$$
In particular, $g$ is $\tau_p$-continuous.
\end{lemma}

\begin{proof}
If $p=\infty$, the result follows from
\cite[Lemma~3]{AzouziDhifaoui2025}, since $\tau_\infty=\bwstar$ and
$\tau_\infty$ is generated by the support functions of $\infty$-compact disks.

Assume now that $1\leq p<\infty$. By sufficiently rapid norm approximation, write
$g=\sum_{n=1}^{\infty}v_n$ with $v_n\in\FVL{E}$ and
$\sum_n2^n\|v_n\|_{\FBLp{p}{E}}<\infty$.
Lemma~\ref{lem:majorant-norm} gives finite majorants
$h_{K_n}=(\sum_j|\delta_{x_{n,j}}|^p)^{1/p}$ such that
\[
 |v_n|\leq h_{K_n},\qquad
 \left(\sum_j\|x_{n,j}\|^p\right)^{1/p}
 \leq2\|v_n\|_{\FBLp{p}{E}};
\]
take $h_{K_n}=0$ if $v_n=0$. Thus
$h_M=(\sum_n2^{np}h_{K_n}^p)^{1/p}$ is the support function of a
$p$-compact disk. With $V_N=\sum_{n=1}^Nv_n$, H\"older gives
\[
 |g-V_N|\leq\sum_{n>N}h_{K_n}\leq2^{-N}h_M,
\]
because the $\ell_{p'}$-norm of $(2^{-n})_{n>N}$ is at most $2^{-N}$. Choose $N$ so large that
$2^{-N}(2h_M(x_0^*)+1)<\varepsilon/2$.
Since $V_N$ is $\tau_p$-continuous, there is a $p$-compact disk
$L\subset E$ such that
\[
 h_L(x^*-x_0^*)<1
 \quad\Longrightarrow\quad
 |V_N(x^*)-V_N(x_0^*)|<\varepsilon/2.
\]
Choose $K$ with $h_K=(h_M^p+h_L^p)^{1/p}$. If
$h_K(x^*-x_0^*)<1$, then $h_M(x^*-x_0^*)<1$ and
$h_L(x^*-x_0^*)<1$, so $h_M(x^*)\leq h_M(x_0^*)+1$. Therefore
\[
 \begin{aligned}
 |g(x^*)-g(x_0^*)|
 &\leq2^{-N}\bigl(h_M(x^*)+h_M(x_0^*)\bigr)
      +|V_N(x^*)-V_N(x_0^*)|\\
 &\leq2^{-N}\bigl(2h_M(x_0^*)+1\bigr)
      +\varepsilon/2
 <\varepsilon.
 \end{aligned}
\]
\end{proof}

For $1\leq p<\infty$, the approximation estimate used in the proof of
Lemma~\ref{lem:summable-local-control} also yields a characterization of
membership in $\FBLp{p}{E}$.

\begin{corollary}\label{cor:summable-approximation}
\leanref{Support}{1902}{summable_approximation}
Let $g\in\Hph{E}$ and $1\leq p<\infty$. Then $g\in\FBLp{p}{E}$ if and only if
there are $V_N\in\FVL{E}$, a $p$-compact disk $K\subset E$, and positive
numbers $\varepsilon_N\to0$ such that
\[
 |g-V_N|\leq \varepsilon_N h_K
 \qquad(N\in\mathbb N).
\]
Moreover, in the forward implication one may take
$\varepsilon_N=2^{-N}$.
\end{corollary}

\begin{proof}
The forward implication is contained in the proof of
Lemma~\ref{lem:summable-local-control}: with the notation used there,
$V_N=\sum_{n=1}^Nv_n$ and the support function
$h_M=(\sum_n2^{np}h_{K_n}^p)^{1/p}$ satisfy
\[
 |g-V_N|\leq2^{-N}h_M.
\]
Conversely, since $h_K\in\FBLp{p}{E}$, monotonicity of the norm
gives
\[
 \|g-V_N\|_{\FBLp{p}{E}}
 \leq \varepsilon_N\|h_K\|_{\FBLp{p}{E}}\longrightarrow0.
\]
Thus $g$ belongs to the norm closure of $\FVL{E}$ in $H_p[E]$, namely
$\FBLp{p}{E}$.
\end{proof}

\begin{corollary}\label{cor:summable-global-domination}
\leanref{Support}{1960}{summable_global_domination}
Let $1\leq p\leq\infty$. For every $g\in\FBLp{p}{E}$ there is a
$p$-compact disk $K\subset E$ such that $|g|\leq h_K$.
\end{corollary}

\begin{proof}
Apply Lemma~\ref{lem:summable-local-control} at zero with
$\varepsilon=1$. If $t>h_K(x^*)$, then
$h_K(x^*/t)<1$, so $|g(x^*)|<t$ by positive homogeneity.
Letting $t\downarrow h_K(x^*)$ gives the assertion, including when
$h_K(x^*)=0$.
\end{proof}

The following provides a convenient characterization of which decreasing nets have infimum 0 in $\FBLp{p}{E}$ in terms of support functions of $p$-compact disks.

\begin{proposition}\label{prop:positive-patches}
\leanref{Support}{2090}{positive_patches}
Let $1\leq p\leq\infty$ and let
$\mathcal A\subset\FBLp{p}{E}_+$ be nonempty. The following are equivalent:
\begin{enumerate}
\renewcommand{\labelenumi}{(\roman{enumi})}
\item there is $0<h\in\FBLp{p}{E}$ such that $h\leq f$ for every
$f\in\mathcal A$;
\item there are $x_0^*\in E^*$, $\varepsilon>0$, and a $p$-compact disk
$K\subset E$ such that
\[
 h_K(x^*-x_0^*)<1
 \quad\Longrightarrow\quad
 f(x^*)\geq\varepsilon
 \qquad(f\in\mathcal A).
\]
\end{enumerate}
Consequently, a decreasing net $(f_\alpha)\subset\FBLp{p}{E}_+$ satisfies
$f_\alpha\downarrow0$ if and only if no such $x_0^*$, $\varepsilon$, and $K$
exist simultaneously for all $\alpha$.
\end{proposition}

\begin{proof}
Assume first that $0<h\leq f$ for every $f\in\mathcal A$. Choose
$x_0^*\in E^*$ with $h(x_0^*)>0$ and put
$\varepsilon=h(x_0^*)/2$. Lemma~\ref{lem:summable-local-control} gives a
$p$-compact disk $K$ such that
\[
 h_K(x^*-x_0^*)<1
 \quad\Longrightarrow\quad
 |h(x^*)-h(x_0^*)|<\varepsilon.
\]
Hence $f(x^*)\geq h(x^*)>\varepsilon$ for every $f\in\mathcal A$ whenever
$h_K(x^*-x_0^*)<1$.

Conversely, suppose that (ii) holds. Taking $x^*=x_0^*$ in (ii) gives
$f(x_0^*)\geq\varepsilon$ for every $f\in\mathcal A$, and hence
$x_0^*\neq0$. Choose $e\in E$ with $x_0^*(e)=1$ and let
\[
 Qx=x-x_0^*(x)e,\qquad M=Q(K).
\]
The disk $M$ is again $p$-compact, and
\[
 h_M(x^*)=h_K\bigl(x^*-x^*(e)x_0^*\bigr).
\]
Set
\[
 h=\varepsilon(\delta_e-h_M)^+\in\FBLp{p}{E}_+.
\]
Since $h_M(x_0^*)=0$, we have $h(x_0^*)=\varepsilon$, so $h\neq0$. If
$h(x^*)>0$, put $a=x^*(e)$. Then $a>h_M(x^*)\geq0$, and therefore
\[
 h_K\left(\frac{x^*}{a}-x_0^*\right)
 =\frac{h_M(x^*)}{a}<1.
\]
Thus $f(x^*/a)\geq\varepsilon$ for every $f\in\mathcal A$. Positive
homogeneity gives
\[
 f(x^*)\geq a\varepsilon\geq h(x^*)
 \qquad(f\in\mathcal A).
\]
The same inequality is trivial when $h(x^*)=0$, proving (i).

For a decreasing positive net, having infimum zero is equivalent to having no
nonzero positive common lower bound, so the final assertion follows from the
equivalence just proved.
\end{proof}

\begin{lemma}\label{lem:subsets-summable-disk}
\leanref{Support}{2378}{subsets_summable_disk}
Let $1\leq p<\infty$. If $A$ is a nonempty subset of a $p$-compact disk
$K\subset E$, then $h_A$ belongs to the norm-closed sublattice of
$\FBLp{p}{E}$ generated by $\{\delta_a:a\in A\}$. More precisely, there
are increasing nonempty finite sets $A_N\subset A$ such that
\[
 s_N=\bigvee_{a\in A_N}|\delta_a|\uparrow h_A
 \quad\text{and}\qquad
 \|h_A-s_N\|_{\FBLp{p}{E}}\longrightarrow0.
\]
\end{lemma}

\begin{proof}
Write $K$ as in \eqref{eq:summable-disk}. Fix $N\in\mathbb N$ and
$\eta>0$. Partition $[-1,1]^N$ into finitely many sets of diameter at most
$\eta$ for the maximum norm. From each cell met by the first $N$
coordinates of a representation of an element of $A$, choose one such
representation, and denote the corresponding elements of $A$ by
$a_1,\ldots,a_m$.

Every $a\in A$ has a representation whose first $N$ coefficients differ
by at most $\eta$ from those of one of these representatives. The
difference of their remaining coefficient sequences has $\ell_{p'}$-norm
at most $2$. Consequently, H\"older gives
\[
 0\leq h_A-\bigvee_{r=1}^m|\delta_{a_r}|
 \leq\eta\sum_{j\leq N}|\delta_{x_j}|
       +2\left(\sum_{j>N}|\delta_{x_j}|^p\right)^{1/p}.
\]
The right-hand side has norm at most
\[
 \eta\sum_{j\leq N}\|x_j\|
 +2\left(\sum_{j>N}\|x_j\|^p\right)^{1/p}.
\]
Since $0\leq h_A\leq h_K$, these estimates are valid in the part of
$\Hph{E}$ with finite $p$-norm. Letting first $N$ increase and then
$\eta$ decrease gives finite supports with arbitrarily small error.
Taking successive unions of their
generating sets yields an increasing sequence $(s_N)$ with
\[
 \|h_A-s_N\|_{\FBLp{p}{E}}\longrightarrow0.
\]
By the monotone convergence lemma, $s_N\uparrow h_A$.
\end{proof}

\subsection{Boundedness below implies order continuity}

Given an embedding $\iota:F\rightarrow E$, we first record the precise closure statement needed for the adjoint
restriction map. Note that this is analogous to the argument based on
Dieudonn\'e's theorem in \cite{AzouziDhifaoui2025}.

\begin{lemma}\label{lem:restriction-polar}
\leanref{Support}{2562}{restriction_polar}
Let $K\subset E$ be a nonempty, absolutely convex, norm-compact set, and put
$D=\{x^*\in E^*:h_K(x^*)<1\}$.  Then
\begin{equation}\label{eq:restriction-polar}
 \overline{\iota^*(D)}^{\,\|\cdot\|_{F^*}}
 =\{y^*\in F^*:h_{K\cap F}(y^*)\leq1\}.
\end{equation}
\end{lemma}

\begin{proof}
Take absolute polars in the dual pairs $(E^*,E^{**})$ and
$(F^*,F^{**})$, identifying $E$ and $F$ with their canonical images.
The set $K$ is weak-star compact in $E^{**}$, so the bipolar theorem gives
$D^\circ=K$.  It follows that
\[
 (\iota^*D)^\circ=(\iota^{**})^{-1}(K)=K\cap F.
\]
For the last equality, use
$\iota^{**}(F^{**})\cap E=F$: an element of this intersection annihilates
$F^\perp$, and $F$ is closed.  Since $\iota^*D$ is absolutely convex, its
bipolar is its norm closure in $F^*$.  This gives
\eqref{eq:restriction-polar}.
\end{proof}

Now, we can prove one implication of our main theorem.

\begin{proposition}\label{prop:poe-sufficiency}
\leanref{Extension}{258}{poe_sufficiency}
If $\overline\iota\colon\FBLp{p}{F}\to\FBLp{p}{E}$ is bounded below, then
it is order continuous.
\end{proposition}

\begin{proof}
Suppose that $f_\alpha\downarrow0$ in $\FBLp{p}{F}$ and
$0<g\leq\overline\iota f_\alpha$ for all $\alpha$.  Choose $x_0^*\in E^*$
with $g(x_0^*)>0$, put $y_0^*=\iota^*x_0^*$, and let
$\varepsilon=g(x_0^*)/2$.  Necessarily $y_0^*\neq0$.  By
Lemma~\ref{lem:summable-local-control}, there is a $p$-compact disk $K\subset E$
such that, with $D=\{x^*:h_K(x^*)<1\}$,
\[
 g(x_0^*+x^*)\geq\varepsilon\qquad(x^*\in D).
\]
Hence
\begin{equation}\label{eq:positive-restriction}
 f_\alpha(y_0^*+\iota^*x^*)\geq\varepsilon
 \qquad(x^*\in D,\ \alpha).
\end{equation}

By Lemma~\ref{lem:subsets-summable-disk}, the support function of $K\cap F$,
viewed on $E^*$, belongs to the norm closure of
$\overline\iota(\FVL{F})$.  Boundedness below makes the range of
$\overline\iota$ closed.  Therefore
\[
 \ell(y^*)=h_{K\cap F}(y^*)\qquad(y^*\in F^*)
\]
defines an element of $\FBLp{p}{F}_+$.  Every $f_\alpha$ is norm continuous on
$F^*$, since it is weak-star continuous on bounded sets.  Thus
Lemma~\ref{lem:restriction-polar} and \eqref{eq:positive-restriction} imply
\begin{equation}\label{eq:positive-source-neighbourhood}
 f_\alpha(y_0^*+y^*)\geq\varepsilon
 \qquad\bigl(\ell(y^*)\leq1,\ \alpha\bigr).
\end{equation}

By Corollary~\ref{cor:summable-global-domination}, choose a $p$-compact disk
$L\subset F$ such that $\ell\leq h_L$. Hence
\[
 h_L(y^*-y_0^*)<1
 \quad\Longrightarrow\quad
 f_\alpha(y^*)\geq\varepsilon
 \qquad(\alpha).
\]
Proposition~\ref{prop:positive-patches} now gives a nonzero positive common
lower bound for the net $(f_\alpha)$, contradicting
$f_\alpha\downarrow0$.
\end{proof}

\subsection{Order continuity forces boundedness below}

To prove the converse implication, we first show that a failure of the norm inequality
can be witnessed by seminorms, rather than arbitrary lattice combinations.
The finite-dimensional intersection in the following proof is key.
For $p=1$ it consists of a polytope.  For $p>1$ it need not be a polytope, but
an arbitrarily close outer polytope approximation gives the required
lattice combination without losing the norm estimate.

\begin{lemma}\label{lem:finite-support-majorant}
\leanref{Support}{2724}{finite_support_majorant}
For every $v\in\FVL{F}$ and every $\eta>0$, there are $z_1,\ldots,z_m\in F$
such that
\[
 r=\bigvee_{k=1}^m|\delta_{z_k}|\geq|v|,
 \qquad
 \|\overline\iota r\|_{\FBLp{p}{E}}
 <\|\overline\iota v\|_{\FBLp{p}{E}}+\eta.
\]
\end{lemma}

\begin{proof}
The assertion is immediate for $v=0$, so assume $v\neq0$.  By
Lemma~\ref{lem:majorant-norm}, choose $x_1,\ldots,x_s\in E$ such that
\[
 |\overline\iota v|\leq
 \left(\sum_{j=1}^s|\delta_{x_j}|^p\right)^{1/p},
 \qquad
 A:=\left(\sum_{j=1}^s\|x_j\|^p\right)^{1/p}
 <\|\overline\iota v\|_{\FBLp{p}{E}}+\frac\eta2.
\]
Put
\[
 K=\left\{\sum_{j=1}^s t_jx_j:
 \|(t_j)_{j=1}^s\|_{\ell_{p'}^s}\leq1\right\},\qquad L=K\cap F.
\]
Both sets are absolutely convex and norm compact in finite-dimensional
subspaces.  We claim that
\begin{equation}\label{eq:finite-support-extension}
 h_L(y^*)=\inf\{h_K(x^*):x^*\in E^*,\ \iota^*x^*=y^*\}
 \qquad(y^*\in F^*).
\end{equation}

Every extension gives the inequality from left to right.  For the reverse
inequality, let $D=\operatorname{span}K$ and equip $D$ with the norm whose
unit ball is $K$.  Hahn--Banach extends $y^*|_{F\cap D}$ to a functional
$w^*\in D^*$ with $h_K(w^*)=h_L(y^*)$.  The functionals $y^*$ and $w^*$
agree on $F\cap D$ and define a continuous functional on $F+D$.  To check
continuity, choose a complement $D_0$ of $F\cap D$ in $D$; then
$F+D=F\oplus D_0$ is a topological direct sum, since $F$ is closed and
$D_0$ is finite-dimensional.  A further Hahn--Banach extension to $E$
proves \eqref{eq:finite-support-extension}.

Since $\overline\iota v=v\circ\iota^*$, the majorant inequality and
\eqref{eq:finite-support-extension} give $|v|\leq h_L$, so $L\neq\{0\}$.
Choose $\theta>0$ such that
$(1+\theta)A<\|\overline\iota v\|_{\FBLp{p}{E}}+\eta$.  In the finite-dimensional
space $\operatorname{span}L$, choose a symmetric polytope $L_\theta$ with
\[
 L\subseteq L_\theta\subseteq(1+\theta)L.
\]
Such outer approximations exist for every finite-dimensional absolutely
convex compact body: take a sufficiently fine finite net in its unit
ball, form the symmetric convex hull, and dilate it slightly.
If $z_1,\ldots,z_m\in F$ are its vertices, then
\[
 r=h_{L_\theta}=\bigvee_{k=1}^m|\delta_{z_k}|\in\FVL{F}_+,
 \qquad |v|\leq h_L\leq r.
\]
On $E^*$, we have
$\overline\iota r\leq(1+\theta)h_K$ because
$L_\theta\subseteq(1+\theta)K$.  Hence $p$-convexity gives
\[
 \|\overline\iota r\|_{\FBLp{p}{E}}
 \leq(1+\theta)\|h_K\|_{\FBLp{p}{E}}
 \leq(1+\theta)A<\|\overline\iota v\|_{\FBLp{p}{E}}+\eta.
\]
\end{proof}

The upper-bound criterion of \cite[Proposition~5.8]{IaninaOikhberg2026}
gives an order-nullness test for the compact support functions used
below. As before, $\|h_K\|_{\FBLp{p}{E}}$ is defined by \eqref{eq:fblp-norm} even
when $h_K$ does not belong to $\FBLp{p}{E}$.

\begin{lemma}\label{lem:compact-support-order-null}
\leanref{Support}{2899}{compact_support_order_null}
Let $E$ be a Banach space, $1\leq p<\infty$, and let $K\subset E$ be
nonempty, absolutely convex, and norm compact.  Suppose that
$h_K$ satisfies $\|h_K\|_{\FBLp{p}{E}}=\infty$ and that
$q_N\in\FBLp{p}{E}_+$ increases pointwise to $h_K$.  Then, for every
$e\in E$,
\[
 (|\delta_e|-q_N)^+\downarrow0\quad\text{in }\FBLp{p}{E}.
\]
\end{lemma}

\begin{proof}
Choose a countable dense subset $D$ of $K$ and define
$T\colon\ell_1(D)\to E$ by $Ta=\sum_{d\in D}a_d d$.  Since $K$ is
absolutely convex and closed, $T(B_{\ell_1(D)})\subseteq K$; hence $T$ is
compact.  Moreover,
\[
 \|T^*x^*\|_{\ell_\infty(D)}
 =\sup_{d\in D}|x^*(d)|=h_K(x^*)\qquad(x^*\in E^*).
\]
By \eqref{eq:fblp-norm}, the absolutely $p$-summing norm of $T^*$ is
$\|h_K\|_{\FBLp{p}{E}}=\infty$.  The restriction of $T^*$ to any
finite-codimensional subspace of $E^*$ also fails to be $p$-summing:
otherwise a finite-rank complement would make $T^*$ $p$-summing.  In the
terminology of \cite[Section~5]{IaninaOikhberg2026}, $T$ is dual
$p$-persistent relative to $E$.

For $e\neq0$ put $u=|\delta_e|$ and $\phi_N=u\wedge q_N$.  If
$g\in\FBLp{p}{E}$ is an upper bound of the $\phi_N$, then pointwise
$g\geq u\wedge h_K$.  Applying
\cite[Proposition~5.8]{IaninaOikhberg2026} with
$x_0=e/\|e\|$ and $\mu=\|e\|$ yields $g\geq u$.  Consequently
$\bigvee_N\phi_N=u$ in $\FBLp{p}{E}$, and
$(u-q_N)^+=u-\phi_N\downarrow0$.  The case $e=0$ is immediate.
\end{proof}

\begin{proposition}\label{prop:poe-necessity}
\leanref{Extension}{429}{poe_necessity}
If $\overline\iota\colon\FBLp{p}{F}\to\FBLp{p}{E}$ is not bounded below,
there are $f_N\in\FVL{F}_+$ with $f_N\downarrow0$ in $\FBLp{p}{F}$ such that
$\overline\iota f_N$ converges in norm to a nonzero positive element of
$\FBLp{p}{E}$.  In particular, $\overline\iota$ is not
$\sigma$-order continuous.
\end{proposition}

\begin{proof}
By density and failure of boundedness below, choose $v_n\in\FVL{F}$ with
$\|v_n\|_{\FBLp{p}{F}}\geq n$ and
$\|\overline\iota v_n\|_{\FBLp{p}{E}}<2^{-n-1}$.  Applying
Lemma~\ref{lem:finite-support-majorant} with $\eta=2^{-n-1}$ gives
seminorms $r_n\in\FVL{F}_+$, each a finite supremum of moduli of
evaluations, such that
\begin{equation}\label{eq:poe-norm-separation}
 \|r_n\|_{\FBLp{p}{F}}\geq n,\qquad
 \|\overline\iota r_n\|_{\FBLp{p}{E}}\leq2^{-n}.
\end{equation}
The series $Q=\sum_{n=1}^{\infty}\overline\iota r_n$ converges in
$\FBLp{p}{E}$, with $\|Q\|_{\FBLp{p}{E}}\leq1$.  Hahn--Banach and
\eqref{eq:poe-norm-separation} give
\[
 \|r_n\|_{\FBLp{\infty}{F}}
 =\|\overline\iota r_n\|_{\FBLp{\infty}{E}}\leq2^{-n}.
\]
Write $r_n=\bigvee_{j=1}^{m_n}|\delta_{z_{n,j}}|$ and let
$K_n=\operatorname{absconv}\{z_{n,1},\ldots,z_{n,m_n}\}$.  Then
$r_n=h_{K_n}$ and $K_n\subseteq2^{-n}B_F$.  The set
\[
 K=\left\{\sum_{n=1}^{\infty}x_n:x_n\in K_n\text{ for all }n\right\}
\]
is absolutely convex and norm compact: the series converge uniformly
over $\prod_nK_n$, and each $K_n$ is compact.  Hence
$h_K =\sum_nh_{K_n}=\sum_nr_n$ pointwise on $F^*$, while
$q_N=\sum_{n=1}^Nr_n\in\FVL{F}_+$ increases to $h_K$.
Since $h_K\geq r_n$ for every $n$, \eqref{eq:poe-norm-separation}
implies $\|h_K\|_{\FBLp{p}{F}}=\infty$.

Choose $e\in F$ with $\|e\|>1$ and set
$f_N=(|\delta_e|-q_N)^+\in\FVL{F}_+$.  By
Lemma~\ref{lem:compact-support-order-null}, $f_N\downarrow0$ in
$\FBLp{p}{F}$.  On the other hand,
$\overline\iota q_N\to Q$ in $\FBLp{p}{E}$, so continuity of lattice
operations gives
\[
 \overline\iota f_N\longrightarrow
 G=(|\delta_e|-Q)^+\quad\text{in }\FBLp{p}{E}.
\]
Here $\|G\|_{\FBLp{p}{E}}\geq\|e\|-\|Q\|_{\FBLp{p}{E}}>0$, because
$|\delta_e|=G+(|\delta_e|\wedge Q)$ and
$0\leq|\delta_e|\wedge Q\leq Q$.  Also
$0<G\leq\overline\iota f_N$ for every $N$.
\end{proof}

\begin{proof}[Proof of Theorem~\ref{thm:poe-inclusion}]
The implication (i)$\Rightarrow$(ii) is immediate,
(ii)$\Rightarrow$(iii) is Proposition~\ref{prop:poe-necessity}, and
(iii)$\Rightarrow$(i) is Proposition~\ref{prop:poe-sufficiency}.
The equivalence of (iii) and (iv) follows from
\cite[Theorem~3.7]{OikhbergEtAl2026}.  The final assertion is the explicit
construction in Proposition~\ref{prop:poe-necessity}.
\end{proof}

\begin{proof}[Proof of Theorem~\ref{thm:intro-poe}]
The assertion is immediate for $T=0$.  If $T\neq0$ and $\overline T$ is
order continuous, Proposition~\ref{prop:necessary} shows that $T$ is an
embedding.  Factor $T=\iota\circ T_0$, where
$T_0\colon F\to T(F)$ is an isomorphism and
$\iota\colon T(F)\hookrightarrow E$ is the inclusion.  The extension
$\overline{T_0}$ is a lattice isomorphism.  Hence $\overline T$ is order
continuous if and only if $\overline\iota$ is order continuous, and
Theorem~\ref{thm:poe-inclusion} gives precisely POE-$p$ for $(T(F),E)$.
The same factorization proves sufficiency.

For the equivalent formulation, boundedness below of $\overline T$ implies
boundedness below of $T$ by evaluation on generators:
\[
 \|Tx\|=\|\overline T\delta_x\|_{\FBLp{p}{E}}\geq c\|\delta_x\|_{\FBLp{p}{F}}
 =c\|x\|\qquad(x\in F)
\]
for some $c>0$.  Apply the factorization above and the equivalence of
boundedness below and order continuity for $\overline\iota$.
Finally, $\overline T=0$ if and only if $T=0$.
\end{proof}

\subsection{Consequences}

The characterization applies, in particular, to every locally complemented
subspace, since these pairs have POE-$p$ for every finite $p$
\cite[Corollary~4.13]{OikhbergEtAl2026}.  Together with
Theorem~\ref{thm:intro-infty}, it therefore extends the
complemented-subspace argument of \cite[Theorem~3.4]{OikhbergEtAl2026}
to locally complemented subspaces for all $1\leq p\leq\infty$.
Likewise, an arbitrary embedding of a $C(K)$ space induces an
order-continuous homomorphism between the corresponding free $p$-convex
Banach lattices for every $p$, by
\cite[Corollary~4.12]{OikhbergEtAl2026} and Theorem~\ref{thm:intro-infty}.

\begin{corollary}\label{cor:bidual-order-continuity}
\leanref{Extension}{992}{bidual_order_continuity}
For every Banach space $E$ and every $1\leq p\leq\infty$, the canonical embedding
$\jmath_E\colon E\to E^{**}$ induces an order-continuous isometric lattice
embedding
\[
 \overline{\jmath_E}\colon\FBLp{p}{E}\longrightarrow\FBLp{p}{E^{**}}.
\]
\end{corollary}

\begin{proof}
Every operator $S\colon E\to\ell_p^n$ extends to
$S^{**}\colon E^{**}\to\ell_p^n$ with the same norm, using the canonical
identification $(\ell_p^n)^{**}=\ell_p^n$.  Thus the pair $(E,E^{**})$ has
POE-$p$ with constant $1$; see also
\cite[Proposition~4.10]{OikhbergEtAl2026}.  The isometric assertion follows
from \cite[Theorem~3.7]{OikhbergEtAl2026}, and order continuity follows
from Theorem~\ref{thm:poe-inclusion} for finite $p$ and
Theorem~\ref{thm:intro-infty} for $p=\infty$.
\end{proof}

It is pertinent to clarify that although as noted in \cite{IaninaOikhberg2026}, the argument of \cite[Theorem 3.6]{AzouziBenRjebTradacete2025} for the non-separable setting uses order continuity of the embedding $\overline{\iota}:FBL[F]\rightarrow FBL[E]$ when $\iota:F\rightarrow E$ denotes an embedding of a separable ideal, the above justifies that this step is now flawless.

At the endpoint, every closed pair $(F,E)$ has POE-$\infty$ with constant
$1$: extend the coordinates of an operator into $\ell_\infty^n$ by
Hahn--Banach.  Thus Theorems~\ref{thm:intro-infty} and
\ref{thm:intro-poe} together give the same POE-$p$ characterization for
all $1\leq p\leq\infty$; only at $p=\infty$ is the extension condition
automatic.


\section{The canonical inclusion between different convexity parameters}

\label{sec:canonical-inclusion}

Let $1\leq p<q\leq\infty$.  Since $\FBLp{q}{E}$ is $p$-convex with
constant one, the universal property yields a contractive lattice
homomorphism
\[
 J_E^{p,q}\colon\FBLp{p}{E}\longrightarrow\FBLp{q}{E},
 \qquad J_E^{p,q}\delta_x=\delta_x\quad(x\in E).
\]
In the functional representations this map is the identity, so it is
injective and has norm-dense range.  In particular,
\[
 \|v\|_{\FBLp{\infty}{E}}\leq\|v\|_{\FBLp{q}{E}}\leq\|v\|_{\FBLp{p}{E}}
 \qquad(v\in\FVL{E}).
\]

 Corollary~\ref{cor:summable-global-domination} connects
the order properties of $J_E^{p,q}$ with containment of compact disks (see also \cite[Proposition 9.38]{OikhbergEtAl2026}).

\begin{theorem}\label{thm:cross-parameter-inclusion}
\leanref{Parameters}{157}{cross_parameter_inclusion}
Let $E$ be a Banach space and let $1\leq p<q\leq\infty$.
The following are equivalent:
\begin{enumerate}
\renewcommand{\labelenumi}{(\roman{enumi})}
\item $J_E^{p,q}$ is order continuous;
\item $J_E^{p,q}$ is $\sigma$-order continuous;
\item $J_E^{p,q}$ is onto;
\item the norms $\|\cdot\|_{\FBLp{p}{E}}$ and
$\|\cdot\|_{\FBLp{q}{E}}$ are equivalent on $\FVL{E}$;
\item $\FBLp{p}{E}$ is $q$-convex;
\item the range of $J_E^{p,q}$ is majorizing in $\FBLp{q}{E}$;
\item the range of $J_E^{p,q}$ is order dense in $\FBLp{q}{E}$;
\item every $q$-compact disk in $E$ is contained in a $p$-compact disk.
\end{enumerate}
\end{theorem}

For the proof we will use separation of the corresponding norms with support functions as follows.

\begin{proposition}\label{prop:cross-parameter-witness}
\leanref{Parameters}{78}{cross_parameter_witness}
Let $E$ be a Banach space and let $1\leq p<q\leq\infty$.
If $J_E^{p,q}$ is not onto, there is a nonempty, absolutely convex,
norm-compact set $K\subset E$ whose support function $h_K$ satisfies
\[
 h_K\in\FBLp{q}{E}_+,\qquad
 \|h_K\|_{\FBLp{q}{E}}\leq1,\qquad \|h_K\|_{\FBLp{p}{E}}=\infty.
\]
Moreover, there are $f_N\in\FVL{E}_+$ and $0<G\in\FBLp{q}{E}$ such that
\[
 f_N\downarrow0\quad\text{in }\FBLp{p}{E},\qquad
 J_E^{p,q}f_N\longrightarrow G\quad\text{in }\FBLp{q}{E},
 \qquad G\leq J_E^{p,q}f_N\quad(N\in\mathbb N).
\]
\end{proposition}

\begin{proof}
Since $J_E^{p,q}$ is injective with dense range, failure of surjectivity
implies failure of boundedness below.  By density, choose
$v_n\in\FVL{E}$ with
\[
 \|v_n\|_{\FBLp{p}{E}}=n,\qquad \|v_n\|_{\FBLp{q}{E}}<2^{-n-1}.
\]

Suppose first that $q<\infty$.  Lemma~\ref{lem:majorant-norm}, at
parameter $q$, gives finite collections $(x_{n,j})_{j=1}^{m_n}$ such that
\[
 |v_n|\leq h_{K_n}=
 \left(\sum_{j=1}^{m_n}|\delta_{x_{n,j}}|^q\right)^{1/q},
 \qquad
 \left(\sum_{j=1}^{m_n}\|x_{n,j}\|^q\right)^{1/q}<2^{-n}.
\]
Concatenating these finite collections gives a $q$-compact disk $K$ with
\[
 h_K=\left(\sum_{n=1}^\infty h_{K_n}^q\right)^{1/q},
 \qquad
 \|h_K\|_{\FBLp{q}{E}}\leq
 \left(\sum_{n=1}^\infty2^{-nq}\right)^{1/q}\leq1.
\]
By Lemma~\ref{lem:subsets-summable-disk}, applied at parameter $q$,
there are finite suprema $s_N\in\FVL{E}_+$ such that
\[
 s_N\uparrow h_K\quad\text{and}\qquad
 \|h_K-s_N\|_{\FBLp{q}{E}}\longrightarrow0.
\]

Suppose now that $q=\infty$, put $c_n=\|v_n\|_{\FBLp{\infty}{E}}$ and choose a
finite-dimensional subspace $M_n\subset E$ containing the generators
of $v_n$.  Hahn--Banach, applied to $x^*|_{M_n}$, gives
\[
 |v_n(x^*)|\leq c_n\|x^*|_{M_n}\|
 =h_{K_n}(x^*),\qquad K_n=c_nB_{M_n}.
\]
Indeed, $v_n$ depends only on this restriction, which extends to $E$
with the same norm.  The set
\[
 K=\left\{\sum_{n=1}^\infty z_n:z_n\in K_n\text{ for all }n\right\}
\]
is absolutely convex and norm compact: each $K_n$ is compact and
the series converge uniformly over $\prod_nK_n$, since
$\sum_nc_n<1$.  Its support function is
\[
 h_K=\sum_{n=1}^\infty h_{K_n},
 \qquad \|h_K\|_{\FBLp{\infty}{E}}\leq\sum_nc_n<1.
\]

For each $N$, choose a finite $2^{-N}$-net $B_N\subset K$, put
$A_N=B_1\cup\cdots\cup B_N$, and set
$s_N=\bigvee_{x\in A_N}|\delta_x|.$
Then
\[
 \|h_K-s_N\|_{\FBLp{\infty}{E}}
 =
 \sup_{\|x^*\|\leq1}|h_K(x^*)-s_N(x^*)|
 \leq 2^{-N}\longrightarrow0.
\]
Since $(s_N)$ is increasing, it follows that
$s_N\uparrow h_K$.

In both cases $h_K\geq|v_n|$ for every $n$, so monotonicity of the
 norm gives $\|h_K\|_{\FBLp{p}{E}}=\infty$.  Choose $e\in E$ with
$\|e\|=2$ and put
\[
 f_N=(|\delta_e|-s_N)^+\in\FVL{E}_+.
\]
Lemma~\ref{lem:compact-support-order-null}, at parameter $p$,
gives $f_N\downarrow0$ in $\FBLp{p}{E}$.  In the target,
\[
 J_E^{p,q}f_N\longrightarrow G=(|\delta_e|-h_K)^+.
\]
Since $|\delta_e|=G+(|\delta_e|\wedge h_K)$, we have
$\|G\|_{\FBLp{q}{E}}\geq2-\|h_K\|_{\FBLp{q}{E}}\geq1$.
Also $s_N\leq h_K$ gives $G\leq J_E^{p,q}f_N$ for every $N$.
\end{proof}

\begin{proof}[Proof of Theorem~\ref{thm:cross-parameter-inclusion}]
The implications (iii)$\Rightarrow$(i)$\Rightarrow$(ii) are immediate.
Conditions (iii) and (iv) are equivalent because $J_E^{p,q}$ has dense
range, and boundedness below makes its range closed.  Condition (iii)
also implies (v), since $q$-convexity is preserved under lattice
isomorphisms.  Conversely, if (v) holds, the universal property gives a
bounded lattice homomorphism
$S\colon\FBLp{q}{E}\to\FBLp{p}{E}$ fixing every $\delta_x$.
Both $SJ_E^{p,q}$ and $J_E^{p,q}S$ agree with the identity on the respective
dense copies of $\FVL{E}$, so $S$ is the inverse of $J_E^{p,q}$.
The implications (iii)$\Rightarrow$(vi) and (iii)$\Rightarrow$(vii)
follow because under (iii) the range is all of $\FBLp{q}{E}$.

Suppose that (iii) fails.  Proposition~\ref{prop:cross-parameter-witness}
gives $h_K\in\FBLp{q}{E}_+$ with $\|h_K\|_{\FBLp{p}{E}}=\infty$,
a sequence $f_N\downarrow0$ in $\FBLp{p}{E}$ and
$0<G\in\FBLp{q}{E}$ with $G\leq J_E^{p,q}f_N$ for every $N$.
The images do not decrease to zero, so (ii) fails.
If (vi) held, then $h_K\leq J_E^{p,q}u$ for some
$u\in\FBLp{p}{E}_+$, contradicting $\|h_K\|_{\FBLp{p}{E}}=\infty$.
Thus (vi) fails.
If (vii) held, order density would give $0<h\in\FBLp{p}{E}$ with
$J_E^{p,q}h\leq G$, and hence $0<h\leq f_N$ for every $N$,
contradicting $f_N\downarrow0$.  Thus (vii) fails.

For (iii)$\Rightarrow$(viii), let $K$ be a $q$-compact disk. Since
$h_K\in\FBLp{q}{E}$, (iii) gives
$h_K\in\FBLp{p}{E}$.
Corollary~\ref{cor:summable-global-domination}, at parameter $p$,
gives a $p$-compact disk $L$ with $h_K\leq h_L$.
Hahn--Banach separation yields $K\subseteq L$.

For (viii)$\Rightarrow$(vi), let $g\in\FBLp{q}{E}$.
Corollary~\ref{cor:summable-global-domination}, at parameter $q$,
gives a $q$-compact disk $K$ with $|g|\leq h_K$.
By (viii), $K\subseteq L$ for some $p$-compact disk $L$.
Hence $g\leq|g|\leq h_K\leq h_L$, with $h_L\in\FBLp{p}{E}$,
which is exactly (vi).
\end{proof}


When $q>2$, the above only holds in trivial cases:

\begin{corollary}\label{cor:cross-parameter-above-two}
\leanref{Parameters}{449}{cross_parameter_above_two}
Let $E$ be a Banach space and let $1\leq p<q\leq\infty$ with $q>2$.
The equivalent conditions of Theorem~\ref{thm:cross-parameter-inclusion}
hold if and only if $E$ is finite-dimensional.
\end{corollary}

\begin{proof}
For nonzero finite-dimensional $E$, choose $x_1,\ldots,x_m\in B_E$ and
$C>0$ such that \[
\|x^*\|
\leq
C\max_{1\leq i\leq m}|x^*(x_i)|
\qquad (x^*\in E^*).
\]
For every family $(x_j^*)_{j=1}^n\subset E^*$ admissible in \eqref{eq:fblp-norm}
with $p=1$, we then have
\[
\sum_{j=1}^n \|x_j^*\|
\leq
C\sum_{i=1}^m\sum_{j=1}^n |x_j^*(x_i)|
\leq
Cm.
\]
Consequently,
\begin{equation}\label{eq:finite-dimensional-norm-comparison}
 \|f\|_{\FBLp{1}{E}}
 \leq
 Cm\|f\|_{\FBLp{\infty}{E}}
 \qquad(f\in\FVL{E}).
\end{equation}
The canonical norm inequalities give equivalence of all the free norms.
The case $E=\{0\}$ is immediate.

If $E$ is infinite-dimensional, \cite[Proposition~9.30]{OikhbergEtAl2026}
shows that $q$-convexity of $\FBLp{p}{E}$ requires
$q\leq\max\{2,p\}$.  This is excluded by the hypotheses.
\end{proof}

\begin{remark}\label{rem:cross-parameter-hilbert}
\leanref{Parameters}{484}{cross_parameter_hilbert}
The restriction $q>2$ in Corollary~\ref{cor:cross-parameter-above-two}
cannot be omitted.  By \cite[Corollary~9.47]{OikhbergEtAl2026}, for
$E=\ell_2$ and $1\leq p<q\leq2$, the canonical map $J_E^{p,q}$ is a
lattice isomorphism onto, hence order continuous.  Thus, in this range
of parameters, infinite dimensionality alone does not imply failure of
order continuity; the criterion is the additional $q$-convexity of
$\FBLp{p}{E}$ in Theorem~\ref{thm:cross-parameter-inclusion}.
\end{remark}

\section*{Funding}
Research partially supported by grants PID2024-162214NB-I00 and CEX2023-001347-S (funded by MCIN/AEI/10.13039/501100011033).

\section*{Acknowledgements}

We thank Timur Oikhberg for discussions that helped us improve the presentation of this paper, and in particular for sharing an argument that simplified an earlier version of Lemma \ref{lem:compact-support-order-null}.

\bibliographystyle{amsplain}
\bibliography{refs}

@book{AliprantisBurkinshaw2006,
  author    = {C. D. Aliprantis and O. Burkinshaw},
  title     = {Positive Operators},
  publisher = {Springer},
  address   = {Dordrecht},
  year      = {2006}
}

@book{AlbiacKalton2016,
  author    = {F. Albiac and N. J. Kalton},
  title     = {Topics in Banach Space Theory},
  edition   = {Second},
  series    = {Graduate Texts in Mathematics},
  volume    = {233},
  publisher = {Springer},
  year      = {2016}
}

@article{AvilesTradacete2023,
  author    = {Avil{\'e}s, A. and Tradacete, P.},
  title     = {Amalgamation and Injectivity in {B}anach Lattices},
  journal   = {Int. Math. Res. Not. (IMRN)},
  volume    = {2023},
  number    = {2},
  pages     = {956--997},
  year      = {2023},
  doi       = {10.1093/imrn/rnab285}
}

@article{ART18,
  author    = {Avil{\'e}s, A. and Rodr{\'i}guez, J. and Tradacete, P.},
  title     = {The free {Banach} lattice generated by a {Banach} space},
  journal   = {J. Funct. Anal.},
  year      = {2018},
  volume    = {274},
  number    = {10},
  pages     = {2955--2977},
  doi       = {10.1016/j.jfa.2018.03.001}
}

@article{AvilesPlebanekRodriguez2018,
  author  = {A. Avil{\'e}s and G. Plebanek and J. D. Rodr{\'i}guez Abell{\'a}n},
  title   = {Chain conditions in free {B}anach lattices},
  journal = {J. Math. Anal. Appl.},
  volume  = {465},
  number  = {2},
  pages   = {1223--1229},
  year    = {2018},
  doi     = {10.1016/j.jmaa.2018.05.061}
}

@article{AzouziDhifaoui2025,
  author        = {Y. Azouzi and W. Dhifaoui},
  title         = {Order denseness in free {B}anach lattices},
  journal       = {Proc. Amer. Math. Soc.},
  year          = {in press},
  eprint        = {2508.11648},
  archiveprefix = {arXiv},
  primaryclass  = {math.FA},
  doi           = {10.1090/proc/17895}
}

@article{AzouziBenRjebTradacete2025,
 author = {Azouzi, Youssef and Rjeb, Asma Ben and Tradacete, Pedro},
 title = {The strong {Nakano} property in {Banach} lattices},
 fjournal = {Revista de la Real Academia de Ciencias Exactas, F{\'{\i}}sicas y Naturales. Serie A: Matem{\'a}ticas. RACSAM},
 journal = {Rev. R. Acad. Cienc. Exactas F{\'{\i}}s. Nat., Ser. A Mat., RACSAM},
 issn = {1578-7303},
 volume = {119},
 number = {4},
 pages = {19},
 note = {Id/No 99},
 year = {2025},
 language = {English},
 doi = {10.1007/s13398-025-01764-7},
 zbMATH = {8095034},
 Zbl = {1579.46012}
}

@article{baker,
  author = {Baker, K. A.},
  title = {Free Vector Lattices},
  journal = {Canad. J. Math.},
  volume = {20},
  year = {1969},
  pages = {58--66}
}

@article{GLTT,
  author = {García-Sánchez, E. and Leung, D. and Taylor, M. and Tradacete, P.},
  title = {{B}anach Lattices with Upper {$p$}-Estimates: Free and Injective Objects},
  journal = {Math. Ann.},
  volume = {391},
  year = {2025},
  pages = {3363--3398}
}

@article{dHT,
  author = {de Hevia, D. and Tradacete, P.},
  title = {Free Complex {B}anach Lattices},
  journal = {J. Funct. Anal.},
  volume = {284},
  number = {10},
  year = {2023},
  pages = {109888}
}

@article{HeinrichMankiewicz1982,
  author  = {S. Heinrich and P. Mankiewicz},
  title   = {Applications of ultrapowers to the uniform and {L}ipschitz classification of {B}anach spaces},
  journal = {Studia Math.},
  volume  = {73},
  pages   = {225--251},
  year    = {1982}
}

@misc{IaninaOikhberg2026,
  author = {A. Ianina and T. Oikhberg},
  title = {Upper bound properties of free and related {Banach} lattices via operators},
  year = {2026},
  eprint = {2609.11928},
  archivePrefix = {arXiv},
  primaryClass = {math.FA},
  url = {https://arxiv.org/abs/2609.11928},
  note = {arXiv:2609.11928}
}

@article{JardonSanchezEtAl2022,
  author  = {H. Jard{\'o}n-S{\'a}nchez and N. J. Laustsen and M. A. Taylor and P. Tradacete and V. G. Troitsky},
  title   = {Free {B}anach lattices under convexity conditions},
  journal = {Rev. R. Acad. Cienc. Exactas F{\'i}s. Nat. Ser. A Mat. RACSAM},
  volume  = {116},
  number  = {1},
  pages   = {Paper No. 15, 25 pp.},
  year    = {2022},
  doi     = {10.1007/s13398-021-01155-8}
}

@article{LaustsenTradacete2026,
 author = {Laustsen, N. J. and Tradacete, P.},
 title = {Banach lattices of positively homogeneous functions induced by a {Banach} space},
 fjournal = {Proceedings of the Edinburgh Mathematical Society. Series II},
 journal = {Proc. Edinb. Math. Soc., II. Ser.},
 issn = {0013-0915},
 volume = {69},
 number = {2},
 pages = {544--580},
 year = {2026},
 language = {English},
 doi = {10.1017/S0013091525101181},
 zbMATH = {8187658}
}

@book{LuxemburgZaanen1971,
  author    = {W. A. J. Luxemburg and A. C. Zaanen},
  title     = {Riesz Spaces I},
  series    = {North-Holland Mathematical Library},
  volume    = {1},
  publisher = {North-Holland},
  address   = {Amsterdam},
  year      = {1971}
}

@misc{MunozLahozTradacete2026,
 author = {D. Mu{\~n}oz-Lahoz and P. Tradacete},
 title = {Free {Banach} {$f$}-algebras},
 year = {2026},
 howpublished = {Preprint, {arXiv}:2511.13299 [math.{FA}] (2026)},
 url = {https://arxiv.org/abs/2511.13299},
 arXiv = {arXiv:2511.13299}
}

@article{Oikhberg2024,
 author = {Oikhberg, T.},
 title = {Geometry of unit balls of free {Banach} lattices, and its applications},
 fjournal = {Journal of Functional Analysis},
 journal = {J. Funct. Anal.},
 issn = {0022-1236},
 volume = {286},
 number = {8},
 pages = {27},
 note = {Id/No 110351},
 year = {2024},
 language = {English},
 doi = {10.1016/j.jfa.2024.110351},
 zbMATH = {7813598},
 Zbl = {1546.46016}
}

@article{OikhbergEtAl2026,
  author  = {T. Oikhberg and M. A. Taylor and P. Tradacete and V. G. Troitsky},
  title   = {Free {B}anach lattices},
  journal = {J. Eur. Math. Soc. (JEMS)},
  volume  = {28},
  number  = {10},
  pages   = {4387--4514},
  year    = {2026},
  doi     = {10.4171/JEMS/1552}
}

@article{OikhbergEtAl2026Corrigendum,
 author = {Oikhberg, T. and Taylor, M. A. and Tradacete, P. and Troitsky, V. G.},
 title = {Corrigendum to: ``{Free} {Banach} lattices''},
 fjournal = {Journal of the European Mathematical Society (JEMS)},
 journal = {J. Eur. Math. Soc. (JEMS)},
 issn = {1435-9855},
 volume = {28},
 number = {10},
 pages = {4515--4518},
 year = {2026},
 language = {English},
 doi = {10.4171/jems/1802},
 zbMATH = {8221494}
}

@article{dePagterWickstead2015,
 author = {de Pagter, B. and Wickstead, A. W.},
 title = {Free and projective {Banach} lattices},
 fjournal = {Proceedings of the Royal Society of Edinburgh. Section A. Mathematics},
 journal = {Proc. R. Soc. Edinb., Sect. A, Math.},
 issn = {0308-2105},
 volume = {145},
 number = {1},
 pages = {105--143},
 year = {2015},
 language = {English},
 doi = {10.1017/S0308210512001709},
 zbMATH = {6438434},
 Zbl = {1325.46020}
}

@misc{SalgueroTradaceteTrejo2025,
 author = {A. Salguero-Alarc{\'o}n and P. Tradacete and N. Trejo-Arroyo},
 title = {Free quasi-{Banach} lattices},
 year = {2025},
 howpublished = {Preprint, {arXiv}:2512.05273 [math.{FA}] (2025)},
 url = {https://arxiv.org/abs/2512.05273},
 arXiv = {arXiv:2512.05273}
}

@article{SinhaKarn2002,
 author = {Sinha, Deba P. and Karn, Anil K.},
 title = {Compact operators whose adjoints factor through subspaces of~{{\(l_p\)}}},
 fjournal = {Studia Mathematica},
 journal = {Stud. Math.},
 issn = {0039-3223},
 volume = {150},
 number = {1},
 pages = {17--33},
 year = {2002},
 language = {English},
 doi = {10.4064/sm150-1-3},
 zbMATH = {1710017},
 Zbl = {1008.46008}
}

@article{Troitsky2019,
 author = {Troitsky, V. G.},
 title = {Simple constructions of {{\(\text{FBL}(A)\)}} and {{\(\text{FBL}[E]\)}}},
 fjournal = {Positivity},
 journal = {Positivity},
 issn = {1385-1292},
 volume = {23},
 number = {5},
 pages = {1173--1178},
 year = {2019},
 language = {English},
 doi = {10.1007/s11117-019-00657-5},
 zbMATH = {7118398},
 Zbl = {1435.46015}
}

\end{document}